\documentclass[12pt]{amsart}
\usepackage{amsfonts, amssymb, amsmath, amsthm, euscript}
\swapnumbers
\theoremstyle{plain}
\newtheorem{thm}{Theorem}[section]
\newtheorem{lem}[thm]{Lemma}
\newtheorem{prop}[thm]{Proposition}

\theoremstyle{definition}
\newtheorem{setup}[thm]{Setup}

\newtheorem*{Ack}{Acknowledgement}

\theoremstyle{remark}

\newtheorem{note}[thm]{}

\def\taubar{\overline{\tau}}
\def\spec{\text{Spec}}

\def\gr{\text{gr}}
\def\id{\text{id}}
\def\a{\mathfrak{a}}
\def\i{\mathfrak{i}}
\def\j{\mathfrak{j}}
\def\oI{\overline{I}}

\def\tildeA{\widetilde{A}}

\def\nxn{n{\times}n}
\def\Olp{{\mathcal O}_{\lambda,p}}
\def\OlpMn{{\mathcal O}_{\lambda,p}(M_n(k))}
\def\hatOlp{\widehat{\mathcal O}_{\lambda,p}}
\def\haty{\widehat{y}}
\def\hattau{\widehat{\tau}}
\def\hatdelta{\widehat{\delta}}
\def\Zp{\mathbb{Z}_p}

\begin{document}

\title[Skew Power Series Rings]{On Prime Ideals of Noetherian \\ 
  Skew Power Series Rings}

\author{Edward S. Letzter}

\address{Department of Mathematics\\
        Temple University\\
        Philadelphia, PA 19122-6094}
      
      \email{letzter@temple.edu }

      \thanks{Research supported in part by grants
        from the National Security Agency.}



      \begin{abstract} We study prime ideals in skew power series
        rings $T:=R[[y;\tau,\delta]]$, for suitably conditioned
        complete right noetherian rings $R$, automorphisms $\tau$ of
        $R$, and $\tau$-derivations $\delta$ of $R$. Such rings were
        introduced by Venjakob, motivated by issues in noncommutative
        Iwasawa theory. Our main results concern ``Cutting Down'' and
        ``Lying Over.''  In particular, assuming that $\tau$
        extends to a compatible automorphsim of $T$, we prove: If $I$
        is an ideal of $R$, then there exists a $\tau$-prime ideal $P$
        of $T$ contracting to $I$ if and only if $I$ is a
        $\tau$-$\delta$-prime ideal of $R$. Consequently, under the
        more specialized assumption that $\delta = \tau - \id$ (a
        basic feature of the Iwasawa-theoretic context), we can conclude:
        If $I$ is an ideal of $R$, then there exists a prime ideal $P$
        of $T$ contracting to $I$ if and only if $I$ is a $\tau$-prime
        ideal of $R$. Our approach depends essentially on two key
        ingredients: First, the algebras considered are zariskian (in
        the sense of Li and Van Oystaeyen), and so the ideals are all
        topologically closed. Second, topological arguments can be
        used to apply previous results of Goodearl and the author on
        skew polynomial rings.
 \end{abstract}

\maketitle


\section{Introduction} 

Given a commutative ring $C$, the following two properties of the
formal power series ring $C[[x]]$ are obvious and fundamental: (1) If
$P$ is a prime ideal of $C[[x]]$ then $P \cap C$ is prime. (2) If $Q$
is a prime ideal of $C$ then there exists a prime ideal of $C[[x]]$
contracting to $Q$. As elementary as these facts are, their analogues
for noncommutative skew power series rings are not immediately
clear. Our aim in this paper is to establish such analogues in the
setting of the skew power series algebras $R[[y;\tau,\delta]]$, for
suitably conditioned complete right noetherian rings $R$,
automorphisms $\tau$ of $R$, and $\tau$-derivations $\delta$ of $R$.
Skew power series rings of this type were introduced by Venjakob in
\cite{Ven}, motivated by issues in noncommutative Iwasawa theory.
Further results were then established in
\cite{SchVen1},\cite{SchVen2}, and several questions on the ideal
theory of these rings were set forth in \cite{ArdBro}.  In \cite{Wan},
completions of certain quantum coordinate rings were also shown to be
iterated skew power series rings in the above sense.

\begin{note} To briefly describe the objects of our study, let $R$ be
  a ring containing an ideal $\i$ such that the $\i$-adic filtration
  of $R$ is separated and complete. Also suppose that the graded ring
  $\gr R$ associated to the $\i$-adic filtration is noetherian; it
  follows that $R$ is noetherian and that the filtration is
  \emph{zariskian\/} in the sense of \cite{LiVOy}.

  Now assume that $\tau$ is an automorphism of $R$, that $\tau(\i)
  \subseteq \i$, that $\delta$ is a left $\tau$-derivation of $R$, and
  that $\tau$ induces a compatible automorphism $\taubar$ of $\gr
  R$. Further assuming that $\delta(R) \subseteq \i$ and $\delta(\i)
  \subseteq \i^2$, we can adapt \cite{SchVen2},\cite{Ven} to construct
  the right noetherian skew power series ring $T:=R[[y;\tau,\delta]]$,
  which will also be complete, separated, and zariskian with respect
  to the $\j$-adic filtration, where $\j = \i + \langle y \rangle$ is
  contained in the Jacobson radical of $T$.  (In this paper we work in
  a slightly different context than in \cite{SchVen2} or
  \cite{Ven}. In \cite{SchVen2} it is assumed that $R$ is not
  necessarily noetherian but is pseudocompact, and in \cite{Ven} it is
  assumed that $R$ is local and $\i$ is the Jacobson radical.)
\end{note}

\begin{note} \label{intro_exampleI} The motivating examples for the
  above and related constructions, in
  \cite{ArdBro},\cite{SchVen1},\cite{SchVen2}, and \cite{Ven}, are the
  noncommutative Iwasawa algebras $\Lambda(G)$, for compact $p$-adic
  Lie groups $G$ containing closed normal subgroups $H$ such that $G/H
  \cong \Zp$. In this case, $\Lambda(G) \cong
  \Lambda(H)[[t;\sigma,\delta]]$, for an automorphism $\sigma$ and
  $\sigma$-derivation $\delta$ of $\Lambda(H)$. We also have
  $\Omega(G) \cong \Omega(H)[[t;\sigma,\delta]]$; again see
  \cite{ArdBro},\cite{SchVen1},\cite{SchVen2}, and \cite{Ven}.
\end{note}

\begin{note} \label{intro_exampleII} Following \cite{Wan}, another
  class of examples arises from the well-known quantized coordinate
  rings that can be realized as iterated skew polynomial rings
\[k[y_1][y_2;\tau_2,\delta_2]\cdots[y_n;\tau_n,\delta_n]\]
over a field $k$, for suitable automorphisms $\tau_i$ and
$\tau_i$-derivations $\delta_i$. (Constructions of these rings can be
found, e.g., in \cite{BroGoo}.)  In \cite{Wan} it is shown, under
addtional assumptions generalizing Venjakob's original approach in
\cite{Ven}, that the completion at $\langle y_1,\ldots,y_n \rangle$ is
an iterated skew power series ring
\[k[[\haty_1]][[\haty_2;\hattau_2,\hatdelta_2]]\cdots[[\haty_n;\hattau_n,\hatdelta_n]],\]
where each $\haty_i$ represents the image of $y_i$ in the completion.
Applicable examples include the quantized coordinate rings of $\nxn$
matrices, of symplectic spaces, and of euclidean spaces; see
\cite{Wan} for further details. \end{note}

\begin{note} To sketch our results, assume for the remainder of this
  introduction that $\tau$ extends to an automorphism of the skew
  polynomial ring $R[y;\tau,\delta]$ and further extends to an
  automorphism of $T = R[[y;\tau,\delta]$. This assumption is
  satisfied by the examples cited above; see (\ref{special_case}).
  Next, recall the definitions of \emph{$\tau$-prime\/} and
  \emph{$\tau$-$\delta$-prime\/} ideals, reviewed in
  (\ref{alpha-prime}) and (\ref{tau-delta-prime}). In particular, the
  presence of notherianity ensures that a $\tau$-prime ideal will be
  the intersection of a finite $\tau$-orbit of prime ideals; see
  (\ref{alpha-prime}ii).
\end{note}

\begin{note} In (\ref{exist}) we prove: (i) {\em If $Q$ is a
    $\tau$-$\delta$-prime ideal of $R$ then there exists a
    $\tau$-prime ideal $P$ of $T$ such that $P \cap R = Q$.} (ii) {\em
    If $P$ is a $\tau$-prime ideal of $T$ then $P\cap R$ is
    a $\tau$-$\delta$-prime ideal of $R$.}  \smallskip
\end{note}

\begin{note} Next, assume further that $\delta = \tau - \id$, a
  condition present for the Iwasawa algebras noted above. Under this
  assumption, an ideal $I$ of $R$ is a $\delta$-ideal if and only if
  $I$ is a $\tau$-ideal, if and only if $I$ is a
  $\tau$-$\delta$-ideal. In (\ref{tau-id-exist}) we prove: (i) {\em If
    $P$ is a prime ideal of $T$ then $P \cap R$ is a $\tau$-prime
    ideal of $R$.} (ii) {\em If $Q$ is a $\tau$-prime ideal of $R$
    then there exists a prime ideal $P$ of $T$ such that $P \cap R =
    Q$.}
\end{note} 

\begin{note} Notation and Terminology. All rings mentioned will be
  assumed to be associative and to posess a multiplicative identity
  $1$. We will use $\spec(\;)$ to denote the prime spectrum of a
  ring. 
\end{note} 

\begin{Ack} My thanks to Konstantin Ardakov, for detailed and helpful remarks on a preliminary draft of this note.
\end{Ack}

\section{Complete Noetherian Rings} \label{complete}

\begin{setup} \label{firstsetup} Throughout, $A$ will denote a
nonzero right noetherian topological ring. We will further assume:

  (1) The topology for $A$ is determined
  by a filtration of ideals
  \[A = \a_0 \supset \a_1 \supset \a_2 \supset \cdots .\]
  (In other words, this filtration forms a fundamental system of open
  neighborhoods of $0$. We use this fundamental system to define
  limits in $A$.)

  (2) Under the given topology, $A$ is Hausdorff (equivalently,
  $\langle 0 \rangle$ is closed in $A$; equivalently, $\cap_i \a_i =
  \langle 0 \rangle$). We say that the filtration is \emph{separated}. 

  (3) $A$ is complete as a topological space (i.e., Cauchy sequences
  converge). We also say that the filtration is complete.

  (4) $A$ contains a dense subring $A'$, equipped
  with the subspace topology.

  (5) $\alpha$ is an automorphism of $A$ restricting to an
  automorphism of $A'$. We include the case when $\alpha$ is the
  identity map.
\end{setup}

When $\a_i = \a^i$ for some ideal $\a$ of $A$, then the preceding
topology on $A$ is the well-known \emph{$\a$-adic topology}. We
include the case when $\a_i = \langle 0 \rangle$ for all $i$, in which
case we obtain the discrete topology on $A$.

\begin{note} \label{alpha-prime} (i) An ideal $I$ of $A$ is an
  \emph{$\alpha$-ideal\/} provided $\alpha(I) \subseteq I$; since $A$
  is right noetherian, note that $\alpha(I) \subseteq I$ if and only
  if $\alpha(I) = I$. We say that $A$ (assumed to be nonzero) is
  \emph{$\alpha$-prime\/} if the product of two nonzero
  $\alpha$-ideals of $A$ is always nonzero, and we say that a proper
  (i.e., not equal to $A$) $\alpha$-ideal $P$ of $A$ is
  \emph{$\alpha$-prime\/} provided $A/P$ is $\alpha$-prime (under the
  induced $\alpha$-action on $A/P$).

  (ii) Since $A$ is right noetherian, an ideal $P$ of $A$ is
  $\alpha$-prime if and only if $P = P_1 \cap \cdots \cap P_t$, for
  some finite $\alpha$-orbit $P_1,\ldots,P_t$, under the naturally
  induced $\alpha$-action on $\spec A$; see \cite[Remarks $4^*$,
  $5^*$, p$.$ 338]{GolMic}. Furthermore, $P_1,\ldots,P_t$ is a
  complete list of the prime ideals of $A$ minimal over $P$.
\end{note}

\begin{note} \label{basic} (See e.g.~\cite[3.2.29]{ArnGlaMik}.) Let
  $B$ denote a complete, Hausdorff topological ring containing a dense
  subring $B'$. If $\varphi: A' \rightarrow B'$ is a bicontinuous ring
  isomorphism, then $\varphi$ extends to a bicontinuous ring
  isomorphism from $A$ onto $B$. In particular, a bicontinuous
  automorphism of $A'$ extends to a bicontinuous automorphism of $A$.
\end{note}

Given a subset $X$ of $A$, we will use $\overline{X}$ to denote the
closure in $A$.

\begin{lem} \label{firstlemma} {\rm (i)} Let $P$ be a closed $\alpha$-prime
  ideal of $A$. Then $P\cap A'$ is an $\alpha$-prime ideal of $A'$.

  {\rm (ii)} If $I$ is an ideal of $A'$ then $\oI$ is an ideal of $A$.

\end{lem}

\begin{proof} (i) Saying that a proper ideal $P$ of $A$ is
  $\alpha$-prime is equivalent to saying, for all $x,y \in A$, and all
  integers $n$, that $\alpha^n(x)Ay \subseteq P$ only if $x$ or $y$ is
  contained in $P$.  Now choose arbitrary $r, s \in A'$, and an
  arbitary integer $n$.  Suppose that $\alpha^n(r)A's \subseteq P\cap
  A'$. To prove $P \cap A'$ is $\alpha$-prime it suffices to prove
  that $r$ or $s$ is contained in $P\cap A'$.

  So let $a$ be an arbitrary element of $A$, and choose a sequence
  $a_i \in A'$ converging to $a$. Then
\[\alpha^n(r)as \; = \; \lim_{i\rightarrow
    \infty} \alpha^n(r)a_is \; \in \; \overline{P\cap A'} \; \subseteq \; P,\]
  because $P$ is closed. Next, since $P$ is $\alpha$-prime and $n$
  was chosen arbitrarily, it follows that one of $r$ or $s$ is
  contained in $P$. Therefore, one of $r$ or $s$ is contained in $P
  \cap A'$, and so $P\cap A'$ is $\alpha$-prime.

(ii) Straightforward.
\end{proof}

\begin{note} \label{zariskianI} We next discuss zariskian filtrations,
  following \cite{LiVOy}.

  (i) To start, we have the \emph{associated graded ring},
\[ \gr A \; := \; A/\a_1 \oplus \a_1/\a_2 \oplus \a_2/\a_3 \oplus \cdots ,\]
and the corresponding function
\[\gr(a) \; := \; (a_i)_{i=0}^\infty, \quad \text{with} \quad a_i \; = \; \begin{cases} a + \a_i/\a_{i+1} &\text{when $a \in \a_i 
    \setminus \a_{i+1}$} \\ 0 &\text{otherwise.}\end{cases} \]
We also have the \emph{Rees ring}:
\[ \tildeA \; := \; A \oplus \a_1 \oplus \a_2 \oplus \cdots \]

(ii) Following \cite[p.~83, Definition]{LiVOy}, we say that $A$ (or
more precisely, the given filtration of $A$) is \emph{right
  zariskian\/} provided that $\a_1$ is contained in the Jacobson
radical of $A$, and provided that $\tildeA$ is right noetherian. Since
the filtration is assumed to be complete, it follows from \cite[p.~87,
Proposition]{LiVOy} that $A$ is right zariskian if $\gr A$ is right
noetherian.

(iii) The consequence of the zariskian property we will require most
is the following: If $A$ is right zariskian then every right ideal of
$A$ is closed with respect to the topology defined by the corresponding
filtration; see \cite[p.~85, Corollary 5]{LiVOy}. 

(iv) Assume that $\gr A$ is right noetherian, and let $I$ be an ideal of $A$. By (iii), the induced filtration on $A/I$,
\[A/I = \a_0 + I \supset \a_1 + I \supset \cdots,\]
is separated and complete. Moreover, $\gr(A/I)$ must be right
noetherian. It follows that $A/I$ is right zariskian with respect to
the induced filtration.

(v) More generally, let $R$ be any ring equipped with a complete
filtration of right ideals. If $\gr R$ is right noetherian then it
follows from well-known arguments that $R$ is also right noetherian.

(vi) Of course, all of the preceding remains true when ``left'' is
substituted for ``right.''
\end{note}

\section{Noetherian Skew Power Series Rings}
\label{SPSR}

We now turn to skew power series rings. Our treatment is adapted from
\cite{SchVen2},\cite{Ven}. However, we do not assume that the
coefficient ring is either local (as in \cite{SchVen2}) or
pseudocompact (as in \cite{Ven}).

\begin{setup} \label{secondsetup} In this section we assume:

  (1) $R$ is a (nonzero) ring containing the proper ideal $\i$.

  (2) The $\i$-adic filtration of $R$ is separated and complete.

  (3) The associated graded ring $\gr R$ corresponding to the
  $\i$-adic filtration is noetherian. Consequently, $R$ is (right and
  left) noetherian, as noted in (\ref{zariskianI}v). Also, by
  (\ref{zariskianI}ii), we know that $R$ is right (and left)
  zariskian, and so all of the right and left ideals of $R$ are closed
  in the $\i$-adic topology, by (\ref{zariskianI}iii).

  (4) $R$ is equipped with an automorphism $\tau$ and a left
  $\tau$-derivation $\delta$ (i.e., $\delta(ab) = \tau(a)\delta(b) +
  \delta(a)b$ for all $a,b \in R$). We assume further that $\tau(\i) =
  \i$ and that $\tau$ induces an automorphism $\taubar$ of $\gr R$,
  with $\taubar \gr (a) = \gr \tau (a)$, for $a \in A$.

  (5) (Following \cite[\S 2]{Ven}.) $\delta(R) \subseteq \i$, and
  $\delta(\i) \subseteq \i^2$. 
\end{setup}

\begin{note} All further references to filtrations and topologies on
  $R$ will refer to the $\i$-adic filtration and $\i$-adic topology.
  Note, since $\tau(\i^i) = \i^i$ for all non-negative integers $i$,
  that $\tau$ is a homeomorphism in the $\i$-adic topology.
\end{note}

\begin{note} \label{skew_poly} We will let $S$ denote the skew
  polynomial ring $R[y;\tau,\delta]$.

    (i) The elements of $S$ are (skew) polynomials
\[ r_0 + r_1 y + r_2 y^2 + \cdots + r_n y^n, \]
for $r_0,\ldots,r_n \in R$, and with multiplication determined by $yr
= \tau(r)y + \delta(r)$, for $r \in R$; more details are given in
(\ref{skew_power_series}ii).

(ii) Set $\tau' = \tau^{-1}$ (recalling that $\tau$ is an
automorphsims), and set $\delta' = -\delta\tau^{-1}$. Then $\delta$ is
a \emph{right\/} $\tau'$-derivation of $R$, and using a process
symmetric to the preceeding we can construct a skew polynomial ring
$S'$ over $R$, with coefficients on the right, such that multiplication is
determined by
\[ ry = y\tau'(r) + \delta'(r) , \]
for $r \in R$. It is well known that $S'$ is isomorphic to $S$.  By
(5) of (\ref{secondsetup}), $\delta'(R) \subseteq \i$ and $\delta'(\i)
\subseteq \i^2$.

(iii) When $\delta$ is the zero derivation, it is customary to denote
the skew polynomial ring by $R[y;\tau]$.

(iv) See, for instance, \cite{GooLet},\cite{GooWar}, or \cite{McCRob}
for detailed background on skew polynomial rings. Recall, in
particular, that $S$ is noetherian since $R$ is noetherian.
\end{note} 

\begin{note} \label{skew_power_series} (Cf.~\cite[\S 1]{SchVen2},
  \cite[\S 2]{Ven}.)

  (i) To start, let $T$ denote the set of (formal skew power) series
\[ \sum r_iy^i \; = \; \sum_{i=0}^\infty r_iy^i ,\]
for $r_0,r_1,\ldots \in R$. 

For now, regard $T$ as a left $R$-module, isomorphic to a direct
product of infinitely many copies of $R$. Identify $S$, as a left
$R$-module, with the left $R$-submodule of $T$ of series having only
finitely many nonzero coefficients.

(ii) For non-negative integers $i$ and all integers $k$, define
$\theta_{ik}(r)$, for $r \in R$, via
\[ y^ir = \theta_{i,i}(r)y^i + \theta_{i, \; i-1}(r)y^{i-1} + \cdots
\theta_{i,0}(r) ,\]
with $\theta_{i,k}(r) = 0$ when $k \leq 0$ or $k > i$. Note that
multiplication in $S$ is given by
\begin{multline}\left( \sum_{i=0}^p a_i y^i \right)\left(\sum_{j=0}^q
    b_j y^j\right) \; = \; \sum_{i=0}^p\sum_{j=0}^q a_i\left(\sum_{k=0}^i\theta_{i,k}(b_j)y^k\right)y^j \; = \\
\sum_{k=0}^\infty \sum_{i=0}^p\sum_{j=0}^q a_i\theta_{i,k}(b_j)y^{k+j} \; = \;
\sum_{n=0}^\infty \sum_{i=0}^p \sum_{j=0}^q a_i\theta_{i,n-j}(b_j)y^n,
\notag\end{multline}
for $a_0,\ldots,a_p,b_0,\ldots,b_q \in R$.

(iii) By (4) and (5) of (\ref{secondsetup}), for $\ell = 2,3,\ldots$, 
\[ \delta(\i^\ell) \; \subseteq \; \i^{\ell -1}\delta(\i) +
\delta(\i^{\ell-1})\i \; \subseteq \; \i^{\ell + 1}.\]
Also, $\theta_{i,k}(R) \subseteq \i^{i-k}$, for all $i \geq k$. Therefore,
\[ \sum_{i=0}^\infty a_i \theta_{ik}(b) \]
converges in $R$, for all $a_0,a_1,\ldots$ and $b$ in $R$.

(iv) By (iii),  we now can (and will) define multiplication in $T$ via
\[ \left( \sum a_i y^ i \right)\left( \sum b_j y^j\right) \; = \;
\sum_{n=0}^\infty\left(\sum_{j=0}^n\left(\sum_{i=0}^\infty a_i\theta_{i,n-j}(b_j)\right)\right)y^n ,\]
for $a_0,b_0,a_1,b_1,\ldots \in R$. It can then be verified that $T$ is an
associative unital ring, and we set $T = R[[y;\tau,\delta]]$. Observe
that $T$ contains $S$ as a subring. (Compare with \cite[\S
1]{SchVen2},\cite[\S 2]{Ven}.)

(v) In view of (\ref{skew_poly}ii), using $\tau'$ and $\delta'$
instead of $\tau$ and $\delta$, we can use a process symmetric to the
preceding to construct an associative skew power series ring $T'$ with
coefficients on the right. In particular, as a right $R$-module, $T'$
is isomorphic to a direct product of infinitely many copies of $R$. It
will follow from (\ref{T=T'}) that $T$ and $T'$ are isomorphic.

(vi) We let $\j$ denote the ideal of $T$ generated by $\i$ and $y$, and we let $\j'$ denote the ideal of $T'$ generated by $\i$ and $y$.
\end{note}

\begin{note} \label{skew_examples} Examples of skew power series rings
  fitting the above framework include:

  (i) (Following \cite[3.3]{ArdBro},\cite[\S 4]{SchVen2},\cite{Ven}.)
  Let $G$ be a compact $p$-adic Lie group containing a closed normal
  subgroup $H$ such that $G/H \cong \Zp$. Then the Iwasawa algebra
  $\Lambda(G)$ has the form $\Lambda(H)[[t;\sigma,\delta]]$, where
  $\sigma$ is an automorphism of the ring $\Lambda(H)$, and where
  $\delta$ is a left $\sigma$-derivation of $\Lambda(H)$. The related
  algebra $\Omega(G)$ can similarly be written as the skew power
  series ring $\Omega(H)[[t;\sigma,\delta]]$. 

  (ii) Following \cite{ArtSchTat} (cf., e.g., \cite{BroGoo} for
  further background and context), let $n$ be a positive integer, let
  $k$ be a field, let $p = (p_{ij})$ be a multiplicatively
  antisymmetric $\nxn$ matrix (i.e., $p_{ij} = p_{ji}^{-1}$ and
  $p_{ii} = 1$) with entries in $k$, and let $\Olp = \OlpMn$ denote
  the multiparameter quantized coordinate ring of $\nxn$
  matrices. Then $\Olp$ is the $k$-algebra generated by
  $y_{11},y_{12},\ldots,y_{nn}$, with the following defining
  relations: $y_{ij}y_{rs} = p_{ir}p_{sj}y_{rs}y_{ij} + (\lambda -
  1)p_{ir}y_{rj}y_{rs}$, for $i > r $ and $j > s$; $y_{ij}y_{rs} =
  \lambda p_{ir}p_{sj}y_{rs}y_{ij}$, for $i > r $ and $j \leq s$;
  $y_{ij}y_{rs} = p_{sj}y_{rs}y_{ij}$, for $i = r $ and $j > s$.

  In \cite{Wan}, the completion $\hatOlp$ of $\Olp$ at $\langle
  y_{11},y_{12},\ldots,y_{nn}\rangle$ is studied, building on the
  approach of \cite{SchVen2},\cite{Ven}.  It is shown that $\hatOlp$
  can be written as an iterated skew power series ring,
  \[k[[\haty_{11}]][[\haty_{12};\hattau_{12},\hatdelta_{12}]] \cdots
  [[\haty_{nn};\hattau_{nn},\hatdelta_{nn}]],\]
  for suitable automorphisms $\hattau_{ij}$ and left
  $\hattau_{ij}$-derivations $\hatdelta_{ij}$. Similar results for
  completions of other quantum coordinate rings can also be found in
  \cite{Wan}.
\end{note}

\begin{lem} \label{jt} For all positive integers $\ell$,
$\j^\ell = \i^\ell + \i^{\ell-1}y + \i^{\ell-2}y^2 + \cdots +
  \i y^{\ell-1} + Ty^\ell$.
\end{lem}

\begin{proof} We prove $\j^\ell \subseteq \i^\ell + \i^{\ell-1}y +
  \cdots + \i y^{\ell-1} + T y^\ell$; the reverse inclusion is
  immediate. To start, recall that $yR \subseteq \i + Ry$. Also, $T =
  R \oplus Ty$. Therefore,
  \[ yT \; = \; y(R \oplus Ty) \; = \; yR + yTy \; \subseteq \; \i +
  Ry + yTy \; \subseteq \; \i + Ty .\]
  So
  \[\j \; = \; \i + \langle y \rangle \; = \; \i + TyT \; \subseteq \i
  + (\i + yT)T \; \subseteq \; \i + \i T + yT \; \subseteq \; \i + yT
  .\]
In particular, the lemma holds for $\ell = 1$. 

Now, for all positive integers $j$, it follows from
(\ref{secondsetup}(5)) that $y \i^j \subseteq \i^jy + \i^{j+1}$.  
Therefore, by induction, for $\ell \geq 2$,
\begin{multline} \j^\ell \; = \; \j^{\ell-1}\j \; \subseteq \;
  (\i + Ty)\left(\i^{\ell-1} + \i^{\ell-2}y + \i^{\ell-3}y^2 + \cdots +
    \i y^{\ell-2} + T y^{\ell-1}\right) \notag \\ \subseteq \;
  \i^\ell + \i^{\ell-1}y + \i^{\ell-2}y^2 + \cdots + \i y^{\ell-1} +
  T y^\ell .
  \end{multline}
The lemma follows.
\end{proof}

\begin{prop} \label{Tcomplete} {\rm (Cf.~\cite[\S 1]{SchVen2},\cite[\S
  2]{Ven}.)} The $\j$-adic filtration on $T$ is separated and
  complete. \end{prop}

\begin{proof} Let $\sum a_i y^i$ be a series in $T$ contained in the
  intersection of all of the powers of $\j$. It then follows from
  (\ref{jt}) that each $a_i$ is contained in the intersection of all
  of the powers of $\i$. Therefore, each $a_i = 0$, since $R$ is
  separated with respect to the $\i$-adic filtration. Hence $T$ is
  separated.

  Now let $\sum a_{i1}y^i, \sum a_{i2}y^i, \ldots$ be a Cauchy
  sequence in $T$ with respect to the $\j$-adic topology. For each
  $i$, it follows that the sequence $\{ a_{ij} \}_{j=0}^{\infty}$ is a
  Cauchy sequence in $R$, by (\ref{jt}), and so we can set
\[ a_i = \lim_{j\rightarrow \infty} a_{ij}, \]
since $R$ is complete with respect to the $\i$-adic topology. Again
using (\ref{jt}), it is not hard to check that
\[ \lim_{j \rightarrow \infty} \sum a_{ij} \; = \; \sum a_i ,\]
under the $\j$-adic topology, and so $T$ is complete with respect to the
$\j$-adic topology.
\end{proof}

\begin{note} \label{S-dense-in-T} Henceforth, references to topologies
  and filtrations on $T$ will refer to the $\j$-adic topology and the
  $\j$-adic filtration. 
By (\ref{jt}), the relative topology on $R$, viewed as a subspace of
  $T$, coincides with the $\i$-adic topology. Note that $S$ is a dense subring of $T$, and equip $S$ with the relative topology from $T$. 
\end{note}

\begin{note} \label{T=T'}
  Consider again $T'$, as described in (\ref{skew_power_series}v), and
  $\j'$ as defined in (\ref{skew_power_series}vi).  By symmetry,
  (\ref{Tcomplete}) tells us that the $\j'$-adic filtration of $T'$ is
  separated and complete, and we henceforth equip $T'$ with the
  $\j'$-adic topology. Next, $S'$, as described in
  (\ref{skew_poly}ii), is a dense subring of $T'$, and we equip $S'$
  with the relative topology. Now observe that the identification of
  $S$ with $S'$ noted in (\ref{skew_poly}iii) is bicontinuous. It then
  follows from (\ref{basic}) and (\ref{S-dense-in-T}) that the ring
  isomorphism from $S$ onto $S'$ noted in (\ref{skew_poly}ii) extends
  to a bicontinuous ring isomorphism from $T$ onto $T'$. We therefore
  identify $T$ with $T'$, and we therefore can write power series in
  $T$ with coefficients on either the left or right.
\end{note}

\begin{note} \label{Tzariskian} Let $\gr T$ denote the associated
  graded ring of $T$ corresponding to the $\j$-adic filtration.

  (i) It follows from (\ref{jt}) that $\gr T \cong (\gr
  R)[y;\taubar]$, where $\taubar$ is the automorphism of $\gr R$ from
  (\ref{secondsetup}(4)). In particular, $\gr T$ is noetherian.

  (ii) Since $\gr T$ is noetherian, and since the $\j$-adic filtration
  of $T$ is complete, by (\ref{Tcomplete}), it follows that the
  $\j$-adic filtration of $T$ is right and left zariskian, by
  (\ref{zariskianI}ii). Consequently the right and the left ideals of
  $T$ are closed in the $\j$-adic topology, by (\ref{zariskianI}iii).
\end{note}

\begin{lem} \label{limit_lemma} Let $\sum a_{i1}y^i, \sum
  a_{i2}y^i, \ldots$ be a convergent sequence in $T$, for $a_{ij} \in
  R$, with
  \[ \lim_{j\rightarrow\infty} \left(\sum_{i=0}^\infty a_{ij}y^i\right) = \sum
  a_iy^i ,\]
  for $a_i \in R$. Then $\lim_{j\rightarrow \infty} a_{ij} = a_i$.
\end{lem}

\begin{proof} Choose a positive integer $\ell$. Then, for some positive
  integer $m$,
\[ \sum (a_i - a_{ij})y^i \; = \; \sum a_iy^i - \sum a_{ij}y^i \; \in \j^\ell,\]
for all $j \geq m$. Next, fix $i$ and choose a positive integer $k$. If
$\ell$ was chosen large enough, it now follows from (\ref{jt}) that
\[ a_i - a_{ij} \in \i^k ,\]
for all $j \geq m$. The lemma follows.
\end{proof}

\begin{prop} \label{cutting_down} {\rm (cf.~\cite[5.12]{GooLet})} Let
  $P$ be a prime ideal of $T$. Then there exists a non-negative
  integer $n$ and a prime ideal $Q$ of $R$ such that $Q,
  \tau(Q),\ldots,\tau^n(Q)$ are exactly the prime ideals of $R$
  minimal over $P\cap R$.
\end{prop}

\begin{proof} It follows from (\ref{firstlemma}i) that $P \cap S$
  is a prime ideal of $S$. Of course, the prime ideals of $R$ minimal
  over $P \cap R$ will be precisely the prime ideals of $R$ minimal
  over $(P\cap S)\cap R$. The proposition now follows from
  \cite[5.12]{GooLet}.
\end{proof}

We now turn to induced ideals.

\begin{note} \label{tau-delta-prime} Recall: (i) Following
  (\ref{alpha-prime}), we have \emph{$\tau$-ideals\/} and
  \emph{$\tau$-prime\/} ideals of $R$.

  (ii) An ideal $I$ of $R$ is a \emph{$\tau$-$\delta$-ideal\/} if
  $\tau(I) = I$ and $\delta(I) \subseteq I$.

  (iii) If the product of nonzero $\tau$-$\delta$-ideals of $R$ is
  always nonzero, we say that $R$ is \emph{$\tau$-$\delta$-prime\/}.
  We say that a proper $\tau$-$\delta$-ideal $Q$ of $R$ is
  \emph{$\tau$-$\delta$-prime\/} provided $R/Q$ is
  $\tau$-$\delta$-prime (under the induced actions of $\tau$ and
  $\delta$ on $R/Q$). Observe that a $\delta$-stable $\tau$-prime
  ideal must be $\tau$-$\delta$-prime.
\end{note}

The following lemma is similar to \cite[2.6--2.8]{LetWan}; however,
the proof is somewhat different.

\begin{lem} \label{induced_ideal} Let $I$ be a $\tau$-$\delta$-ideal of
  $R$. Then: 

{\rm (i)}
  \[I[[y;\tau,\delta]] \; := \; \left\{ \sum a_iy^i \; : \; a_i \in I
  \right\} \; = \left\{ \sum y^ib_i \; : \; b_i \in I
  \right\}.\]

{\rm (ii)} $I[[y;\tau,\delta]]$ is is an ideal of $T$ and is the topological
closure in $T$ of $IS=SI$.

{\rm (iii)} $IT = I[[y;\tau,\delta]] = TI$.

{\rm (iv)} $T/IT \cong (R/I)[[y;\tau,\delta]]$, under the induced
actions of $\tau$ and $\delta$ on $R/I$.

\end{lem}

\begin{proof} (i) Recall $\theta_{ik}$ from
  (\ref{skew_power_series}ii). Since $I$ is $\tau$-$\delta$-stable, it
  follows that $\theta_{ik}(I) \subseteq I$ for all non-negative
  integers $i$ and $k$. Next, for $b_i\in I$, it follows from
  (\ref{skew_power_series}iii) that
\[ \sum_{i=0}^\infty \theta_{ik}(b_i) \]
is well defined for all $k$; the closure of $I$ in $R$, established in
(3) of (\ref{secondsetup}), further ensures that this series is contained
in $I$. Therefore,
\begin{multline} \sum y^i b_i \; = \;
\left(\sum_i \theta_{i0}(b_i)\right) +
  \left(\sum_i \theta_{i1}(b_i)\right)y + \left(\sum_i
    \theta_{i2}(b_i)\right)y^2 + \cdots \notag \\ \in \; \left\{ \sum b'_k y^k
    \; : \; b'_k \in I \right\}.
\end{multline}
A symmetric argument, beginning with $\sum a_iy^i$, for $a_i \in I$,
establishes (i); see (\ref{skew_power_series}v).

(ii) Since $I$ is a $\tau$-$\delta$-ideal and $\tau$ is an automorphism,
\[ IS = \{ a_0 + a_1y + \cdots + a_my^m : a_i \in I\} \; = \; \{ b_0 +
yb_1 + \cdots + y^nb_n : b_i \in I\} \; = \; SI .\]
However, since $I$ is a closed subset of $R$, it follows from
(\ref{limit_lemma}) that $I[[y;\tau,\delta]]$ is the topological
closure in $T$ of the ideal $IS = SI$ of $S$. Furthermore, it now
follows from (\ref{firstlemma}iv) that $I[[y;\tau,\delta]]$ is an ideal
of $T$.

(iii) First, it is easy to see that $IT$ and $TI$ are contained in
$I[[y;\tau,\delta]]$. On the other hand, $IT$ and $TI$ are closed
subsets of $T$, by (\ref{Tzariskian}ii), since $IT$ is a right ideal
of $T$ and $TI$ is a left ideal of $T$. We learned in (ii) that
$I[[y;\tau,\delta]]$ is the topological closure of $IS = SI$, and so
$I[[y;\tau,\delta]]$ is contained in $IT$ and $TI$, since $IS
\subseteq IT$ and $SI \subseteq TI$. Part (iii) follows.

(iv) Straightforward, given (ii) and (iii).
\end{proof}

\begin{note} \label{special_case} A well known special case occurs
  when $\tau$ can be extended to compatible automorphisms of $S$ and
  $T$. (Here, compatible means that $\tau$ extends to an automorphism
  of $T$ that restricts to an automorphism of $S$). This situation
  occurs, for instance, when $\delta\tau = \tau\delta$ (as operators
  on $R$); in this case we can extend $\tau$ to $S$ and $T$ by setting
  $\tau(y) = y$. (A proof of this assertion will follow from the next
  paragraph, setting $q = 1$.) In particular, $\delta$ and $\tau$ will
  satisfy the equation $\delta\tau = \tau\delta$ when $\delta = \tau -
  \id$.  For the examples $\Lambda(H)[[t;\sigma,\delta]]$ and
  $\Omega(H)[[t;\sigma,\delta]]$ mentioned in (\ref{skew_examples}i),
  it is a basic feature of their construction that $\delta = \sigma -
  \id$; see \cite[\S 4]{SchVen2} and \cite[2.2]{Ven}.

  More generally, suppose for the moment that $\delta\tau =
  q\tau\delta$ for some central unit $q$ of $R$ such that $\tau(q) =
  q$ and $\delta(q) = 0$.  Following \cite[2.4ii]{GooLet}, $\tau$
  extends to an automorphism of $S$ such that $\tau(y) = q^{-1}y$. We
  see that this automorphism of $S$ is bicontinuous, and so $\tau$
  extends to an automorphism of $T$, by (\ref{basic}). The condition
  $\delta\tau = q\tau\delta$ holds for the skew power series discussed
  in (\ref{skew_examples}ii).

  When $\tau$ extends to compatible automorphisms of $S$ and $T$, we
  can refer to $\tau$-prime ideals of these rings, following
  (\ref{alpha-prime}).
\end{note}

\begin{note} \label{tau-contract} Assume that $\tau$ extends to
  compatible automorphisms of $S$ and $T$. Let $I$ be a $\tau$-ideal
  of $S$ or $T$. Then, for all $a \in I$, $ya - \tau(a)y \in I$. It
  follows that $I \cap R$ is a $\tau$-$\delta$-ideal of $R$.
\end{note}

Part (i) of the following, a weakened analogue to \cite[3.3i]{GooLet},
establishes ``Lying Over,'' and part (ii) provides another ``Cutting
Down.''

\begin{thm} \label{exist} Assume that $\tau$ extends to compatible
  automorphisms of $S$ and $T$.

  {\rm (i)} Suppose that $Q$ is a $\tau$-$\delta$-prime ideal of $R$.
  Then there exists a $\tau$-prime ideal $P$ of $T$ such that $P \cap
  R = Q$.

  {\rm (ii)} Suppose that $P$ is a $\tau$-prime ideal of $T$. Then
  $P\cap R$ is a $\tau$-$\delta$-prime ideal of $R$. 
\end{thm}

\begin{proof} (i) It follows from (\ref{induced_ideal}) that $QT =
  Q[[y;\tau,\delta]] = TQ$, from which it is easy to deduce that $QT
  \cap R = Q$. Also, $\tau(QT) = \tau(Q)\tau(T) = QT$, and so
  $QT$ is a $\tau$-ideal of $T$. We see, then, that there exists at
  least one $\tau$-ideal of $T$ whose intersection with $R$ is equal
  to $Q$. Therefore, we can choose a $\tau$-ideal $P$ of $T$
  maximal such that $P\cap R = Q$. (Note that $P \ne T$.)

  We can prove as follows that $P$ is $\tau$-prime: Let $I$ and $J$ be
  $\tau$-ideals of $T$, both containing $P$, such that $IJ \subseteq
  P$.  By (\ref{tau-contract}), $I \cap R$ and $J \cap R$ are
  $\tau$-$\delta$-ideals of $R$. Also, $(I \cap R)(J \cap R) \subseteq
  P\cap R = Q$. Now, since $Q$ is $\tau$-$\delta$-prime, at least one
  of $I \cap R$ or $J \cap R$ is contained in $Q$; say $I \cap R
  \subseteq Q$.  But $I \supseteq P$, and so $I \cap R = Q$. The
  maximality of $P$ now ensures that $I = P$.

  We conclude that $P$ is a $\tau$-prime ideal of $T$ contracting to
  $Q$.

(ii) Let $I = P \cap R$. It follows from (\ref{tau-contract}) that $I$ is
a $\tau$-$\delta$-ideal of $R$.

 Now suppose that $J$ and $K$ are $\tau$-$\delta$-ideals of $R$ such
that $JK \subseteq I$. By (\ref{induced_ideal}iii), $TJ = JT$ and $TK
= KT$ are ideals of $T$. Moreover, $JT$ and $KT$ are $\tau$-ideals of
$T$.  Now, $(JT)(KT) = (JK)T \subseteq (P\cap R)T \subseteq P$, and so
one of $JT$ or $KT$ is contained in $P$. By (\ref{induced_ideal}iii),
$JT \cap R = J[[y;\tau,\delta]] \cap R = J$ and $KT \cap R =
K[[y;\tau,\delta]]\cap R = K$. We see that one of $J$ or $K$ is
contained in $I$, and (ii) follows.
\end{proof}

We can obtain even more precise results in the case where $\delta =
\tau - \id$. As noted in (\ref{special_case}), the skew polynomial
rings $\Lambda(H)[[t;\sigma,\delta]]$ and
$\Omega(H)[[t;\sigma,\delta]]$ satisfy $\delta = \sigma - \id$. Note,
when $\delta = \tau - \id$, that an ideal $I$ of $R$ is a
$\delta$-ideal if and only if $I$ is a $\tau$-ideal, if and only if
$I$ is a $\tau$-$\delta$-ideal.

\begin{thm} \label{tau-id-exist} Assume that $\delta = \tau - \id$. 

  {\rm (i)} Let $P$ be a prime ideal of $T$. Then $P \cap R$ is a
  $\tau$-prime ideal of $R$. In particular, $P \cap R = Q \cap \tau(Q)
  \cap \cdots \cap \tau^{n-1}(Q)$, for some prime ideal $Q$ of $R$,
  and some positive integer $n$, such that $\tau^n(Q) = Q$.

  {\rm (ii)} Let $Q$ be a $\tau$-prime ideal of $R$.  Then there
  exists a prime ideal $P$ of $T$ such that $P \cap R = Q$.
\end{thm}

\begin{proof}
Following \cite[4.1]{Ven}, if we set $z = 1+y$, then
\[ zr = (1+y)r = r + yr = r + \tau(r)y + \delta(r) = r +
  \tau(r)y + \tau(r) - r = \tau(r)(1+y) = \tau(r)z,\]
  for $r \in R$. Also, $z=1+y$ is a unit in $T$, with inverse $z^{-1}
  = 1 - y + y^2 - y^3 + \cdots$. We thus obtain an inner automorphism
  of $T$, via $z(\;\;)z^{-1}$. 

  This inner automorphism maps $y$ to itself and restricts to $\tau$
  on $R$. Setting $\tau(\;\;) = z(\; \;)z^{-1}$ on $T$, we obtain
  exactly the extension of $\tau$ to $T$ discussed in
  (\ref{special_case}). Moreover, $\tau:T \rightarrow T$ restricts to
  the unique automorphism of $S$ that maps $y$ to itself and maps each
  $r \in R$ to $\tau(r)$.

  Note, for all ideals $I$ of $T$, that $\tau(I) = zIz^{-1} = I$. In
  other words, every ideal of $T$ is a $\tau$-ideal. Consequently, it
  follows directly from the definition (\ref{tau-delta-prime}ii) that
  an ideal of $T$ is prime if and only if it is $\tau$-prime. Also, by
  (\ref{tau-contract}), every ideal of $T$ contracts to a
  $\tau$-$\delta$-ideal of $R$.

  Next, observe that $S = R[y;\tau,\delta] = R[z;\tau]$. Let $P$ be a
  prime ideal of $T$. We know from (\ref{firstlemma}i) that $P \cap S$
  is a prime ideal of $S$. Also, since $z$ is a unit in $T$, we know
  that $z \notin P \cap S \subseteq P$. It now follows from the
  standard theory that $P \cap R = (P \cap S) \cap R$ is $\tau$-prime;
  see, for example, \cite[10.6.4iii]{McCRob}. Recalling
  (\ref{alpha-prime}ii), part (i) follows.

  Part (ii) now follows from (\ref{exist}i), since an ideal of $R$ is
  $\tau$-prime if and only if it is $\tau$-$\delta$-prime.
\end{proof}

\begin{note} Assume that $\tau$ extends to compatible automorphisms of
  $S$ and
  $T$, and let $Q$ be a $\tau$-$\delta$-prime ideal of $R$. In
  \cite[3.3i]{GooLet} it is proved that $QS = SQ$ is $\tau$-prime.  We
  ask: Is $TQ = QT$ a $\tau$-prime ideal of $T$?
\end{note}



\begin{thebibliography}{99}
  
\bibitem{ArdBro} K. Ardakov, and K. A. Brown, Ring-theoretic
  properties of Iwasawa algebras: a survey, \emph{Doc.~Math.}, Extra
  Vol.~(2006), 7--33.

\bibitem{ArnGlaMik} V. I. Arnautov, S. T. Glavatsky, and
  A. V. Mikhalev, \emph{Introduction to the Theory of Topological
    Rings and Modules}, Pure and Applied Mathematics 197, Marcel
  Dekker, Inc. New York, 1996.

\bibitem{ArtSchTat} M. Artin, W. Schelter, and J. Tate, Quantum
  deformations of $GL_n$, \emph{Communic.~Pure Appl.~Math.}, 44
  (1991), 879-895.

\bibitem{BroGoo} K. A. Brown and K. R. Goodearl, \emph{Lectures on Algebraic
    Quantum Groups}, Advanced Courses in Mathematics CRM Barcelona,
  Birkh\"auser Verlag, Basel, 2002.

\bibitem{GolMic} A. Goldie and G. Michler, Ore extensions and polycyclic
  group rings, \emph{J. London Math.~Soc.~(2)}, 9 (1974/75), 337--345.

\bibitem{Goo} K. R. Goodearl, Prime ideals in skew polynomial rings and
quantized Weyl algebras, \emph{J. Alg$.$}, 150 (1992), 324--377.

\bibitem{GooLet} K. R. Goodearl and E. S. Letzter, Prime ideals in
  skew and $q$-skew polynomial rings, \emph{Mem.~Amer.~Math.~Soc.},
  109 (1994). 

\bibitem{GooWar} K. R. Goodearl and R. B. Warfield, Jr., \emph{An
    Introduction to Noncommutative Noetherian Rings, Second Edition},
  London Mathematical Society Student Texts 61, Cambridge University
  Press, Cambridge, 2004.

\bibitem{LetWan} E. S. Letzter and L. H. Wang, Noetherian Skew Inverse
  Power Series Rings, \emph{Algebr.~Represent.~Theory}, to appear.

\bibitem{LiVOy} H. Li and F. Van Oystaeyen, \emph{Zariskian Filtrations},
  $K$-Monographs in Mathematics 2, Kluwer Academic Publishers, Dordrecht,
  1996.

\bibitem{McCRob} J. C. McConnell and J. C. Robson, \emph{Noncommutative
    Noetherian Rings}, Graduate Studies in Mathematics 30, American
  Mathematical Society, Providence, 2000.

\bibitem{SchVen1} P. Schneider and O. Venjakob, Localisations and
  completions of skew power series rings, arXiv:0711.2669.

\bibitem{SchVen2} \bysame, On the codimension of modules
  over skew power series rings with applications to Iwasawa algebras,
  \emph{J. Pure Appl$.$ Algebra}, 204 (2006), 349--367.

\bibitem{Ven} O. Venjakob, A non-commutative Weierstrass preparation theorem
  and applications to Iwasawa theory (with an appendix by Denis Vogel),
  \emph{J. Reine Angew$.$ Math.}, 559 (2003), 153--191.

\bibitem{Wan} L. Wang, Completions of quantum coordinate rings,
  \emph{Proc.~Amer.~Math.~Soc.}, to appear.

\end{thebibliography}
\end{document}